\documentclass[12pt]{amsart}
\usepackage{latexsym}
\usepackage{psfig}
\usepackage{amsmath,amssymb}

\newtheorem{thm}{Theorem}[section]
\newtheorem{cor}[thm]{Corollary}
\newtheorem{con}[thm]{Conjecture}
\newtheorem{lemma}[thm]{Lemma}
\newtheorem{cl}[thm]{Claim}
\newtheorem{note}[thm]{Note}

\newtheorem{observe}[thm]{Observation}
\newtheorem{prop}[thm]{Proposition}

\numberwithin{equation}{section}

\author{Joseph Cohen}
\title{Primitive Roots in Quadratic Fields}
\email{coheny@techunix.technion.ac.il}
 \address{Math Department,
Technion, Haifa, 32000, Israel}
\date{\today}
\thanks{Submitted in partial fulfillment of the requirements for the degree of Doctor of Philosophy.
 Supported by grants from the Technion-Israel Institute of Technology}
\begin{document}
\maketitle

\begin{abstract}
We consider an analogue of Artin's primitive root conjecture for
units in real quadratic fields.  Given such a nontrivial unit, for
a rational prime $p$ which is inert in the field the maximal order
of the unit modulo $p$ is $p+1$. An extension of Artin's
conjecture is that there are infinitely many such inert primes for
which this order is maximal. This is known at present only under
the Generalized Riemann Hypothesis. Unconditionally, we show that
for any choice of 7 units in different real quadratic fields
satisfying a certain simple restriction, there is at least one of
the units which satisfies the above version of Artin's conjecture.

\end{abstract}

\newcommand{\ord}{\operatorname{ord}}

\section{Introduction}
A natural question to ask is if there are many primes for which
$2$ is a primitive root, that is if the subgroup $\langle 2
\rangle$ of the multiplicative group $\mathbb{F}_p^\star$ of the
field of $p$ elements generated by $2$ is the whole group. Is
there a finite number of such groups $\mathbb{F}_p^\star$?
 Does the same apply for any integer $a$?

In 1927 Emil Artin made the following conjecture (\cite{A}):
\begin{con}
Let $a\neq-1$ be an integer which is not a perfect square. Then
there are infinitely many primes p such that
 $$<a>=\mathbb{F}_p^\star \ .$$
 In addition, for $x>0$, the number of primes $p \leq x$ with this
 property is asymptotic as $x\to \infty$ to
 $$A(a)\frac{x}{\log x}$$
 where $A(a)$ is a constant which depends on $a$.
\end{con}

In 1967 Hooley proved (\cite{H}) Artin`s conjecture with the
asymptotic formula under the Generalized Riemann Hypothesis. In
1983 Gupta and Murty (\cite{GM}) proved that there are 13 specific
integers such that at least one of them fulfills the Artin
conjecture. From the proof we can deduce that Artin`s conjecture
is true for almost all integers. R. Murty, K. Murty and Gupta
(\cite{GMM}) showed that we can reduce the specific set of
integers from 13 to 7. Improving the analytic part of Gupta and
Murty enabled Heath-Brown to give the best result till now:
\begin{thm}\label{theorem 1.2}
Let q, r and s any three primes. Then at least one of them is a
primitive root $mod \ p$ for infinitely many primes p.
\end{thm}

We note that theorem~\ref{theorem 1.2} holds for any three non-zero integers,
$q, r$ and $s$ which are multiplicatively independent  where $q$,
$r$, $s$, $-3qr$, $-3qs$, $-3rs$ and $qrs$ are not a square. (we say that
$r$ integers $a_1,...,a_r$ are multiplicatively independent if for
any integers $n_1,...,n_r, \ a_1^{n_1} \cdot \cdot \cdot
a_r^{n_r}=1 \Rightarrow n_1=...=n_r=0$).

In this work we present an analog of Artin`s conjecture in a
different field and we will prove a result similar to the one just
shown (we will show that a set which contains a specific number of
elements or more always contains a primitive root).

\subsection{Artin`s conjecture in a real quadratic field}
Let $d \neq 1$ be a square-free natural number and let $\Delta=d$
if $d\equiv 1 (mod \ 4)$ and  $\Delta=4d$ otherwise. Let
$K=\mathbb{Q}(\sqrt{\Delta})$ be a real quadratic field and denote
the integer ring of $K$ by $\mathcal{O}_K$.
The principal ideals that are generated by a rational prime $p$,
$p\mathcal{O}_K$ take one of the following forms
\begin{enumerate}
\item $p\mathcal{O}_K=P$  (inert);
\item $p\mathcal{O}_K=P_1P_2, \ \ \ P_1 \neq P_2$ (splits);
\item $p\mathcal{O}_K=P^2$ (ramified)
\end{enumerate}
where $P$ and $P_i$ are prime ideals in $\mathcal{O}_K$.
We note that the option (3) occurs only in a finite number of
cases and so does not interest us. \

Now, the norm map
$$\mathcal{N}:\mathcal{O}_K \mapsto \mathbb{Z}$$
gives a homomorphism
$$\mathcal{O}_K/(p)  \mapsto \mathbb{Z}/p$$
For any unit $\epsilon$ the kernel of this map contains the
residue class $\epsilon$ modulo $p$. Denote this kernel by
$C_\epsilon(p)$. By lemma 19 in \cite{KR2} (appendix B)
$$
\ord (C_\epsilon(p))=\left\{\begin{array}{lr}p-1,& p \ splits\\ p+1,& p \ inert
                    \end{array}\right.
$$

Assuming  GRH,  Cooke and Weinberger (\cite{CW}) and Lenstra
(\cite{Len}) showed that given a real quadratic field $K$, 
there are infinitely many split primes for which the fundamental unit 
of the field\footnote{assume that it has norm $+1$} 
has maximal order (namely $p-1$) in $C_\epsilon(p)$.

Using  the strong analytic theorem of Heath-Brown \cite{H}, 
   Narkiewicz \cite{N2} proved the following unconditional theorem:
 \begin{thm}
 Let $\epsilon_1, \epsilon_2, \epsilon_3$ be units in the
 integers rings $\mathcal{O}_{K_1},\mathcal{O}_{K_2},
 \mathcal{O}_{K_3}$ of $K_1=\mathbb{Q}(\sqrt{\Delta_1})$,
 $K_2=\mathbb{Q}(\sqrt{\Delta_2})$,
 $K_3=\mathbb{Q}(\sqrt{\Delta_3})$, respectively, which are not roots of unity. There is an index j,
$1 \leq j \leq 3$,
 such that for infinitely many split primes $p$,
$c_j\epsilon_j, \ c_j=\pm 1$, has order $p-1$ $(mod \ (p))$.
\end{thm}

For inert primes, one wants similar results. On GRH, an analogue of 
\cite{CW} \cite{Len} was only proven recently by Roskam (\cite{R}). 
We want to extend the result of Narkiewicz for inert primes. 

 In this case the order of $C_\epsilon(p)$ $(mod \
p)$ is $p+1$. So we cannot use the the result of Heath-Brown on
the divisors of $p-1$. We shall use a simpler method to get
infinitely many primes $p$ such that $\frac{p+1}{2}=P_3$ (we write
$P_3$ for an integer with at most three prime factors) but with
almost same magnitude of the prime divisors. With this result we
obtain:

\begin{thm}\label{theorem 1.6}
Let $\epsilon_1,\ldots,\epsilon_7$ be units in the rings of
integers $\mathcal{O}_{\Delta_1},\ldots,\\
\mathcal{O}_{\Delta_7}$ of
${\mathbb{Q}(\sqrt{\Delta_1})},\ldots,{\mathbb{Q}(\sqrt{\Delta_7})}$,
respectively, which are not roots of unity, with
$\Delta_1,\ldots,\Delta_7$ multiplicatively independent and
distinct from $3$. Assume that  all the numbers
$(-1)^{a_1}3^{a_2}\prod\limits_{i=1}^7{\Delta_i^{b_i}}$,   $a_i,
b_i\in \{0, 1\}$, are not perfect squares if
$\sum\limits_{i=1}^7{b_i}$ is odd. Then there exists an index
$1\leq j\leq 7$, such that for infinitely many inert primes $p$,
the unit $c_j\epsilon_j$, ($c_j=\pm 1$), has order $p+1$ modulo
${p}\mathcal{O}_{\Delta_j}$.
\end{thm}

\begin{cor}\label{corollary 1.7}
Let $\epsilon_1,\ldots,\epsilon_7$ be units in the rings of
integers $\mathcal{O}_{\Delta_1},\ldots,\\
\mathcal{O}_{\Delta_7}$ of ${\mathbb{Q}(\sqrt{\Delta_1})},\ldots,
{\mathbb{Q}(\sqrt{\Delta_7})}$,respectively, which are not roots
of unity, with $\Delta_1,\ldots,\Delta_7$ primes distinct from
$3$. Then there exists an index  $1\leq j\leq 7$, such that for
infinitely many inert primes $p$, the unit $ c_j\epsilon_j$,
($c_j=\pm 1$) has order $p+1$ modulo  ${p}
\mathcal{O}_{\Delta_j}$.
\end{cor}

\section{The work of Gupta-Murty and of Heath-Brown}
Since our work is based on the idea of Gupta and Murty with the
advanced version as in the paper of Heath-Brown it will be natural
to present their work. We start with following
trivial idea:
Since the number of elements in $\mathbb{F}_p^\star$ is $p-1$, if
we show for all integer $d$, and infinitely many primes $p$
$$
d\mid p-1 , d \neq p-1 \Rightarrow a^d\neq1 \ \ mod \ p
$$
we will have proven the conjecture.

So our first  goal is to find infinitely many primes $p$ with a
small number of prime divisors of $p-1$. Heath-Brown proved the
following lemma.

\begin{lemma}\label{lemma 1.3}
Let $k=1, 2$ or $3$ and $K=2^k$. Let $q, r$ and $s$ be any three
primes. Then for any sufficiently large $x \in \mathbb{R}^+$ there
exist two numbers $ \ \epsilon, \delta \in (0,1/4)$ and
$c=c(\epsilon,\delta)>0$ so that there are at least
$c\frac{x}{\log^2x}$ primes $p\leq x$ which satisfy:

Either $\frac{p-1}{K}$ is prime or $\frac{p-1}{K}=p_1p_2$ for
$p_1, p_2$ primes $p_i > p^{{1/4}+\epsilon}, \ i=1,2$ and
$p_1<p^{{1/2}-\delta}$. Furthermore, $p$ satisfies
$$(\frac{q}{p})=(\frac{r}{p})=(\frac{s}{p})=-1.$$
\end{lemma}

Now we prove theorem~\ref{theorem 1.2} from this lemma. Assume for simplicity
that $K=2$ and that we have infinitely many primes $p\leq x$ as in
the lemma~\ref{lemma 1.3} such that $\frac{p-1}{2}=l$ where $l$ is a prime.
Take one of the three primes in the lemma, say, $q$. If the order
of $q$ equals $l$ we get a contradiction to the fact that
$(\frac{q}{p})=-1$ by the theorem on cyclic groups (the order of
the squares subgroup of $\mathbb{F}_p^\star$ is
$\frac{p-1}{2}=l$). If $\ord(q)=2l$ we are done. If not, the only
possibility left is $\ord(q)=2$ but this does not occur for
sufficiently large primes $p$ and hence $q$ is a primitive root.

Assume now that there exist $c\frac{x}{\log^2x}$ primes $p\leq x$
as in lemma~\ref{lemma 1.3} such that $\frac{p-1}{2}=p_1p_2$. As before the
order of $q$, $r$ and $s$ can be (if they are not primitive) $2$,
$2p_1$ or $2p_2$. As before there is only a small number of cases
where $\ord(q)=2$. Assume that $\ord(q)=2p_1$. For this case we need
some observation. Let $n$ be a natural number and $\Omega(n)$
denote the number of prime factors of $n$ (with multiplicity) and
write $f \ll g $ (or $g \gg f$ or $f=O(g))$, where $g$ is a
positive function, if there exists a constant $c>0$ such that
$|f(x)|\leq cg(x)$. Then
\begin{observe}\label{observation 1.4}
(1) For any natural number $n$, $\Omega(n)\ll \log n$.

(2) Given an integer $a$, The number of primes p such that $\ord(a)
< y \ (mod \ p)$ is $O(y^2)$.
\end{observe}
To see (1), use
$n=q_1^{\alpha_1} \cdot\cdot\cdot q_r^{\alpha_r}\geq 2^{\alpha_1}
\cdot\cdot\cdot 2^{\alpha_r}\geq
2^{\alpha_1+...+\alpha_r}=2^{\Omega(n)}$.
To see (2), use
$ \sum\limits_{m < y}\Omega(a^{m}-
1)
 \ll_a \sum\limits_{m < y}\log(a^{m}- 1)
 \ll_a \sum\limits_{m < y}m \ll_a y^2 $.

Now, if $\ord(q)=2p_1 < x^{1/2-\delta} \ (mod \ p),$  then by
observation~\ref{observation 1.4} this occurs for at most
$(x^{1/2-\delta})^2=x^{1-2\delta}$ primes.

Since $\frac{x^{1-2\delta}} {cx/\log^2x} \rightarrow 0$ as $x
\rightarrow \infty$ there are a negligible number of primes \ $p$
\ such that $\ord(q)=2p_1 \ (mod \ p).$ This fact is also true
for $r$ and $s$.

Now, assume that $q$ and $r$ and $s$ have order $2p_2$. Since
$\mathbb{F}_p^\star$ cyclic group,
$\ord(<q,r,s>)=2p_2<x^{{3/4}-\epsilon}$. By lemma 2 in \cite{GM}
the number of primes $p$ such that $\ord(<q,r,s>)<y$ is
$O(y^{4/3})$. So the number of primes $p$ such that
$\ord(<q,r,s>)=2p_2<x^{{3/4}-\epsilon}$ is $O(x^{1-{4\epsilon/3}})$
and as before, is negligible in comparison to $\frac{cx}{\log^2x}$.

 \section{Notation and Preliminaries}
Now before we prove the theorem about the prime divisors of $p+1$
(as in lemma~\ref{lemma 1.3} for $p-1$) we need to decide on some notation.

 Let $\Pi(y;m,s)$  denote the number of primes $p \leq x$ such that
 $p\equiv s \ (mod \ m)$ where $m$ and $s$ are some integers, and
$$
E(y;m,s):= \Pi(y;m,s) -\frac{Li{(y)}}{\varphi{(m)}}
$$
where $Li(y)=
 \int_2^y\frac{dt}{\log t}$. Also set
$$
E(x;m):= \max\limits_{1\leq y \leq x} \ \max\limits_{(s,m)=1}|E(y;m,s)| \;.
$$

 Define $\mathcal{A}=\{p+1|p\leq x,p\equiv u \
 (mod \ {v})\}$ where $u,\ v$ are some integers such that $(u,v)=1, \ u\equiv1 (\mod 2),\
8|v$, $(\frac{u+1}{2},v)=1$ and take $X=\frac{Li(x)}{\varphi(v)}$.\\

For a square-free integer $d$,   $(d,v)=1$, let
\begin{equation*}
\begin{split}
\mathcal{A}_d &:=
\{a\in\mathcal{A}:  a\equiv 0  \mod  d)\}\\
&=\{p+1: p\leq x,\ p \equiv u \mod  v,\ p \equiv -1 \mod d \}
\end{split}
\end{equation*}

By the Chinese remainder theorem there exists an $l$ such that
$$
|\mathcal{A}_d| =\#\{p+1|p\leq x,p \equiv l \ (mod \ {dv})\}\;.
$$
By the definition of $E(x;dv,l)$,
$$
|\mathcal{A}_d|=\frac{Lix}{\varphi(dv)}+E(x;dv,l)=
\frac{1}{\varphi(d)}\frac{Lix}{\varphi(v)} +
E(x;dv,l)=\frac{X}{\varphi(d)}+ E(x;dv,l)
$$

Define $\omega(d):=\frac{d}{\varphi(d)}$ and
$$R_d:=|\mathcal{A}_d|-\frac{\omega(d)}{d}X = E(x;dv,l)$$
Finally, we define two arithmetical functions for a square-free
$d=p_1 \cdot\cdot\cdot p_k$. $\mu(d)=(-1)^k$ and $\nu(d)=k$ (where
$\mu(1)=1$ and $\nu(1)=0$).

Now we want to prove two lemmas.

\begin{lemma}\label{lemma 2.1}
for any prime $q$, which is relatively prime to $v$ we have:
\begin{equation}\label{equation 2.1}
 0 \leq \frac{1}{q-1} \leq1-\frac{1}{c_1}
 \end{equation}
 where $c_1>1$ is some suitable constant.

\begin{equation}\label{equation 2.2}
 \sum\limits_{w \leq q < z}
\frac{\log q}{q-1}-\log\frac{z}{w} = O(1) \  \  \ (2 \leq w \leq z)
\end{equation}
where $O$ does not depend on $z$ or $w$.

 \begin{equation}\label{equation 2.3}
\prod\limits_{\substack{2<q<z \\ q\nmid v}}
(1-\frac{1}{q-1})\gg\frac{1}{\log z}.
\end{equation}
 \end{lemma}

\begin{proof}
Since $q>2$ it is clear that \eqref{equation 2.1} holds.

As for the second equation,
 $ \sum\limits_{w\leq q<z}\frac{\log q}{q-1}=\sum\limits_{w\leq q
<z}\frac{\log q}{q}\frac{q}{q-1}=\sum\limits_{w\leq q
<z}\frac{\log q}{q}(1+\frac{1}{q-1})=\sum\limits_{w\leq q
<z}\frac{\log q}{q}+\sum\limits_{w\leq q
<z}\frac{\log q}{q(q-1)}=\log\frac{z}{w}+O(1) \ \ \ \ \
(\sum\limits_{p<x}\frac{\log p}{p}= \log x+O(1))$.

Hence we get \eqref{equation 2.2}. Finally,
\begin{equation*}
\begin{split}
\prod\limits_{\substack{2<q<z \\ q\nmid v}}
(1-\frac{\omega(q)}{q})&=\prod\limits_{\substack{2<q<z \\q \nmid
v}}{(1-\frac{1}{q-1})} \gg
\prod\limits_{2<q<z}{(1-\frac{1}{q-1})} \\
&= \exp({\log{\prod\limits_{2<q<z}{(1-\frac{1}{q-1})}}}) \\
&=
\exp({\sum\limits_{2<q<z}{\log{(1-\frac{1}{q-1})}}}) \\
&\gg
\exp({\sum\limits_{2<q<z}(-\frac{1}{q-1}-\frac{1}{(q-1)^2})})
\end{split}
\end{equation*}

Since
$$
\frac{1}{q-1}=\frac{1}{q}+\frac{1}{q(q-1)} \leq \frac{1}{q} +
\frac{1}{(q-1)^2}
$$
and $\sum\limits_{2<q<z}\frac{1}{(q-1)^2}$ converges, we get
$$
\prod\limits_{2<q<z}(1-\frac{1}{q-1})\gg
\exp({-\sum\limits_{2<q<z}\frac{1}{q}})
$$
Since
$$\sum\limits_{2<q<z}{\frac{1}{q}}\sim \log\log z$$
we have
$$ \exp({-{\sum\limits_{2<q<z}\frac{1}{q}}})\gg
 \exp({-\log\log z})=\frac{1}{\log z}$$
 \end{proof}

\begin{lemma}\label{lemma 2.2}
For any natural square-free number $d$, $(d,v)=1$, given an $A>0$
there exist constants $c_2(\geq 1)$ and $c_3(\geq 1)$ such that
 \begin{equation}\label{equation 2.4}
 \sum\limits_{d<\frac{X^\frac{1}{2}}{(\log x)^{c_2}}}
  \mu^2(d)3^{\nu(d)}|R_d| \leq c_3 \frac{X}{\log^AX}, \ \ \ (X
  \geq2)
\end{equation}
\end{lemma}

\begin{proof}
denote by $S_{R_d}$ the term which we need to estimate.
$$
S_{R_d}=\sum\limits_{d<\frac{X^\frac{1}{2}}{(\log
    x)^{c_2}}}\mu^2(d)3^{\nu(d)}|R_d|
$$
By the definitions of $R_d$ and $E(x;dv)$
$$
S_{R_d}\leq
\sum\limits_{d<\frac{X^\frac{1}{2}}{(\log x)^{c_2}}}
\mu^2(d)3^{\nu(d)}|E(x;dv)|
$$
since $E(x;dv)\ll \frac{x}{dv}$ if $d\leq\frac{x}{v}$ we get
that
$$
S_{R_d} \ll
x^\frac{1}{2}\sum\limits_{{d<\frac{X^\frac{1}{2}}{(\log x)^{c_2}}}}
\frac{\mu^2(d)3^{\nu(d)}}{d^\frac{1}{2}}|E(x;dv)|^{\frac{1}{2}}.
$$
 By Cauchy`s inequality,
$$
S_{R_d}\ll
x^\frac{1}{2}(\sum\limits_{d<X^\frac{1}{2}}\frac{\mu^2(d)3^{2\nu(d)}}{d})^\frac{1}{2}
(\sum\limits_{dv<\frac{vX^\frac{1}{2}}{(\log
    x)^{c_2}}}|E(x;dv)|)^\frac{1}{2}.
$$

For sufficiently large $x$ we obtain
$$
S_{R_d}\ll
x^\frac{1}{2}(\sum\limits_{d<x^\frac{1}{2}}\frac{\mu^2(d)3^{2\nu(d)}}{d})^\frac{1}{2}
(\sum\limits_{dv<\frac{x^\frac{1}{2}}{(\log
    x)^{c_2}}}|E(x;dv)|)^\frac{1}{2}.
$$

 With Bombieri-Vinogradov theorem (\cite{BV})
 (given any positive constant $e_1$, there exist a positive
 constant $e_2$ such that $\sum\limits_{d<\frac{x^\frac{1}{2}}{\log^{e_2}x}}E(x;d)=O(\frac{x}{\log^{e_1}x})$)
for the last sum and since
$\sum\limits_{d<w}\frac{\mu^2(d)9^{\nu(d)}}{d}\leq (\log w
+1)^9$ (see \cite{HR}, p.115, equation (6.7)) we find that for given
 constant $B$ there exist $c_2$ such that
$$S_{R_d}\ll \frac{x}{\log^B{x}}$$
So, for given $A$ there exist $c_2$ such that
$$S_{R_d} \ll \frac{X}{\log^A{X}}$$
where $\ll$ depends on $v$ and $c_2$.
\end{proof}

\section{Proof of the theorem~\ref{theorem 1.6} - the sieve part}
In this section we will show that for a sufficiently small
$0<\delta <1/4$ there exists some constant $c(\delta)>0$ (which
depends on $\delta$) such that for at least $c(\delta)
\frac{x}{\log^2x}$ primes $p\leq x, \ p\equiv u \ (mod \ v), \
\frac{p+1}{2}=P_3$ where $q|\frac{p+1}{2} \Rightarrow q >
x^{1/4-\delta}$. Later we will sharpen this result further.

\subsection{Use of the lower bound linear sieve}
 In the following subsection we will show, using the linear sieve, that for
 a sufficiently small $0 < \delta < 1/4$ there exists some constant $c_1(\delta)>0$
(which depends on $\delta$) such that for at least $c_1(\delta)
\frac{x}{\log^2x}$ primes $p\leq x, \ p\equiv u \ (mod \ v), \
\frac{p+1}{2}$ has at most four prime divisors
  all of them greater than $x^{\frac{1}{4}-\delta}$. \\

Define $S(\mathcal{A},z)=
\#\{a|a\in\mathcal{A},(a,\prod\limits_{\substack{p<z \\
p \nmid v}}p)=1\}$ and let $f$ denote the "lower bound function"
for the linear sieve which is defined as $f(t)=2e^\gamma
t^{-1}\log(t-1)$ for\ $2\leq t \leq4$ , where $\gamma$ is Euler
constant. Then (see \cite[Theorem 8.4, page 236]{HR}):
\begin{lemma}\label{lemma 3.1}
Assume \eqref{equation 2.1}, \eqref{equation 2.2} and \eqref{equation 2.4}.
Then for $X^{1/8}<z<X^{1/4}$ we have
\begin{equation}
 S(\mathcal{A},z)\geq X\prod\limits_{\substack{q<z \\ q\nmid v}} (1-\frac{\omega(q)}{q})
 \{f(\frac{\log x}{2\log z})+O(\frac{1}{\log x})\}
\end{equation}
where the O-term does not depend on $X$ or on $z$.
\end{lemma}

\begin{note} \rm
Obviously $z$ influences the number of primes which divide the
elements of $\mathcal{A}$ and their magnitude. Heath-Brown used a
stronger version of this lemma which gives
$z=x^{{1/4}+\epsilon_0}$ where $\epsilon_0$ is a specific small
real number.

\end{note}

By Lemmas \ref{lemma 2.1} and \ref{lemma 2.2}, \eqref{equation 2.1},
\eqref{equation 2.2} and \eqref{equation 2.4} hold.
Hence we can use lemma~\ref{lemma 3.1} with $z=X^{\frac{1}{4}-\delta}$.
$$
S(\mathcal{A},X^{\frac{1}{4}-\delta})\geq
X\prod\limits_{\substack{q<X^{\frac{1}{4}-\delta}\\ q\nmid v}}
(1-\frac{1}{q-1})\{f(\frac{1}{2}\frac{\log x}{\log x^{\frac{1}{4}-\delta}})
+O(\frac{1}{\log x})\}
$$
by Lemma~\ref{lemma 2.1} \eqref{equation 2.3} we have
$$
S(\mathcal{A},X^{\frac{1}{4}-\delta}) \gg \frac{X}{\log
  x^{\frac{1}{4}-\delta}}f(\frac{2}{1-4\delta})
$$
But for $2\leq t \leq4$, $f(t)=2e^\gamma t^{-1}\log(t-1)$, and so,
\begin{multline*}
S(\mathcal{A},X^{\frac{1}{4}-\delta})\gg
\frac{X}{\log x^{1/4-\delta}}2e^\gamma(\frac{1-4\delta}{2})
\log\frac{1+4\delta}{1-4\delta} \\
\gg \frac{x}{\log^2x}\log
\frac{1+4\delta}{1-4\delta} = \frac{x}{\log^2x}\log
(1+\frac{8\delta}{1-4\delta}).
\end{multline*}

Since $\log(1+s)/s\sim 1$  as  $s\to 0$  and
for, $ 0 < \delta < 1/4$, $1-4\delta$ are bounded, we have:
\begin{lemma}
$$
S(\mathcal{A},X^{\frac{1}{4}-\delta}) \gg \delta
\frac{x}{\log^2x}
$$
where the implied constant in $\gg$ does not depend on $\delta$.
\end{lemma}

\begin{note} \rm By definition of $S(\mathcal{A},X^{\frac{1}{4}-\delta})$, for all
sufficiently small $0<\delta<1/4$, there are at least $\gg \delta
\frac{x}{\log^2x}$ primes $p\leq x$, such that any prime divisor of
$p+1$ (for a $p$ in our sequence) is greater than
$X^{\frac{1}{4}-\delta}$ or divides $v$. Since by our assumption
$(\frac{u+1}{2},v)=1$ where $p \equiv u \ (mod \ {v})$ and
$X=\frac{Li(x)}{\varphi(v)}$ we obtain that all odd prime divisor
of $\frac{p+1}{2}$ are greater than $x^{\frac{1}{4}-\delta}.$
 Hence there are at most four primes divisors of $\frac{p+1}{2}$
 which are greater than $x^{\frac{1}{4}-\delta}$. In the next subsection we will
 show that there is only a small number of primes $p\leq x$  such
 that $\frac{p+1}{2}$ has exactly four primes divisors all of which are greater than $x^{\frac{1}{4}-\delta}.$
\end{note}

\subsection{First use of the Selberg upper bound sieve}

In order to prove that there is only a small number of primes
$p\leq x$ such that exactly four primes divide $\frac{p+1}{2}$ we
need to use Selberg`s upper bound sieve
(see \cite[theorem 3.12]{HR}):
\begin{prop}\label{proposition 3.5}
 Let $a, b$ be integers satisfying
$$
ab\neq 0, \quad \gcd(a,b)=1,\quad  2\mid ab
$$
Then as $x \rightarrow \infty$ we have uniformly in $a, b$ that
\begin{multline*}
|\{p:p \leq x,\ ap+b=\mbox{ prime} \}| \leq \\
8\prod\limits_{p>2}(1-\frac{1}{(p-1)^2})\prod\limits_{2<p|ab}\frac{p-1}{p-2}
         \frac{x}{\log^2x}\{1+O(\frac{\log\log x}{\log x})\}
\end{multline*}
\end{prop}
 From this proposition we derive the following:
 \begin{lemma}\label{lemma 3.6}
For any $0<\delta<1/4$,
there exists  $c_2(\delta)\frac{x}{\log^2x}\ (c_2(\delta) >0)$ primes
$p\leq x$ such that $\frac{p+1}{2}$ has at most three prime divisors
all of  which are greater than $x^{{1/4}-\delta}.$
 \end{lemma}

\begin{proof}
Assume that $\frac{p+1}{2}=p_1p_2p_3p_4,\ p \leq x$ where the
$p_is$ primes greater than $x^{1/4-\delta}$. Instead of counting
the elements in this set we can count the products of primes
$p_1p_2p_3p_4$ such that $2p_1p_2p_3p_4-1=p \leq x$ where the
$p_is$ are primes greater than $x^{1/4-\delta}$.

To count the latter set we use proposition~\ref{proposition 3.5}.
We take $a=2p_1p_2p_3$,  $b=-1$  and $Y=\frac{x+1}{2p_1p_2p_3}$
(since $2p_1p_2p_3p_4-1\leq x  \Leftrightarrow p_4\leq
\frac{x+1}{2p_1p_2p_3}$).

By the proposition,
\begin{multline*}
S_{p_4}=\#\{p_4\leq Y: ap_4+b=\mbox{ prime}\}\\
=  \#\{p_4\leq \frac{x+1}{2p_1p_2p_3}: 2p_1p_2p_3p_4-1=\mbox{ prime}\} \\
\ll\frac{x+1}{2p_1p_2p_3\log^2\frac{x+1}{2p_1p_2p_3}}
\prod\limits_{\substack{p|2p_1p_2p_3\\p\neq2}} \frac{p-1}{p-2}
\end{multline*}

Since the $p_i$'s are big primes,
the term $\prod\limits_{\substack{p|2p_1p_2p_3\\
p\neq2}}\frac{p-1}{p-2}$ is approximately one.
Then
$$S_{p_4} \ll \frac{x+1}{2p_1p_2p_3\log^2\frac{x+1}{2p_1p_2p_3}}$$

 From the fact that for all $i=1,2,3, \ p_i<x^{1/4+3\delta}$ we
have for a sufficiently small $\delta$
\begin{equation*}
  \begin{split}
S_{p_4} & \ll \frac{1}{p_1p_2p_3\log^2\frac{x}{(x^{1/4+3\delta})^3}} \\
& \ll \frac{1}{p_1p_2p_3} \cdot \frac{x}{\log^2{x^{1/4-9\delta}}} \ll
(\frac{1}{1/4-9\delta})^2\frac{1}{p_1p_2p_3} \cdot \frac{x}{\log^2x}
  \end{split}
\end{equation*}

Now we shall sum-up the last term over all possibilities for
$p_1,p_2,p_3$. This number is bounded by
$$S^*_{p_4}=4\frac{x}{\log^2x}\sum\limits_{p_1}\frac{1}{p_1}\sum\limits_{p_2}\frac{1}{p_2}
\sum\limits_{p_3}\frac{1}{p_3}$$
where the sum is over
$x^{1/4-\delta}<p_i<x^{1/4+3\delta}$, $i=1,2,3$.

\begin{observe}\label{observation 3.7}
We have $\sum\limits_{x^\beta<p<x^\alpha}\frac{1}{p}=
\log\frac{\alpha}{\beta}$ +o(1).
\end{observe}

By  observation~\ref{observation 3.7},
$$
S^*_{p_4} \ll
\log^3\frac{1/4+3\delta}{1/4-\delta}\frac{x}{\log^2x}\ll
\log^3(1+\frac{16\delta}{1-4\delta})\frac{x}{\log^2x}
$$

Since $\log(1+s)=O(s)$ for $0<s<1$
and for, $0< \delta < 1/4, 1-4\delta$ is bounded, we have
$$
S^*_{p_4} \ll \delta^3 \frac{x}{\log^2x}
$$
where $\ll$ does not depend on $\delta$.
 Hence, $S^*_{p_4}$ is a small number in comparison to
$S(\mathcal{A},X^{\frac{1}{4}-\delta})\gg
 \delta\frac{x}{\log^2x}$.
\end{proof}

\subsection{Second use of Selberg`s upper bound sieve}
Up till now we know that for any sufficiently small number
$\delta>0$, there are $c_2(\delta)\frac{x}{\log^2x}$ primes $p\leq
x$ such that $\frac{p+1}{2}$ has at most three prime divisors all
of which are greater than $x^{1/4-\delta}$. In this section we
want to prove the existence of $c_3(\delta) \frac{x}{\log^2x}$
primes $p \leq x$, such that $\frac{p+1}{2}=P_3$ and if
$\frac{p+1}{2}$ is a product of exactly three primes $q_3 \geq q_2
\geq q_1$ then $q_1>x^{1/4-\delta}, \ q_2>x^{1/4+2\delta}, \
q_3>x^{1/3+\delta^2}$. First we prove the claim about $q_2$ (by
the previous subsections it is
 clear that $q_1>x^{1/4-\delta}$).

Assume that $q_1$ and $q_2$ take values between
$x^{\frac{1}{4}-\delta}$ and $x^{\frac{1}{4}+2\delta}$. Instead of
count the number of primes $p\leq x$ such that
$\frac{p+1}{2}=q_1q_2q_3$ where $q_1$ and $q_2$ are between
$x^{\frac{1}{4}-\delta}$ and $x^{\frac{1}{4}+2\delta}$ we shall
count the products $q_1q_2q_3$ such that $2q_1q_2q_3-1=p\leq
x$ where $q_1$ and $q_2$ are between $x^{\frac{1}{4}-\delta}$ and
$x^{\frac{1}{4}+2\delta}$.

To count this set we use Proposition~\ref{proposition 3.5}.
Define $a=2q_1q_2, \ b=-1$  and $Y=\frac{x+1}{2q_1q_2}$ (since
$2q_1q_2q_3-1 \leq x$   $\Leftrightarrow q_3\leq
\frac{x+1}{2q_1q_2}$).
By Proposition~\ref{proposition 3.5}
\begin{equation*}
\begin{split}
S_{q_3}& =\#\{q_3\leq Y: aq_3+b=\mbox{ prime}\}\\
& =\#\{q_3\leq \frac{x+1}{2q_1q_2}: 2q_1q_2q_3-1=\mbox{ prime} \}\\
& \ll\frac{x+1}{2q_1q_2\log^2\frac{x+1}{2q_1q_2}}
\prod\limits_{\substack{p\mid 2q_1q_2\\
p\neq 2}}\frac{p-1}{p-2}
\end{split}
\end{equation*}

As in the previous subsection, since the $q_i`s$ are big primes
the term $\prod\limits_{\substack{p|2q_1q_2\\
p\neq2}}\frac{p-1}{p-2}$ is approximately one, so
$$
S_{q_3} \ll \frac{x+1}{2q_1q_2\log^2\frac{x+1}{2q_1q_2}}\;.
$$
Now we sum-up the last term over all possibilities for $q_1,q_2$.
This number bound by, (see the last previous subsection).
$$
S^*_{q_3}= \frac{x}{\log^2x} \sum\limits_{x^{\frac{1}{4}-\delta}
\leq q_1 \leq x^{\frac{1}{4}+2\delta}}\frac{1}{q_1}
\sum\limits_{x^{\frac{1}{4}-\delta} \leq q_2 \leq
x^{\frac{1}{4}+2\delta}}\frac{1}{q_2}\ll\frac{x}{\log^2x}\log^2\frac{1+8\delta}{1-4\delta}
$$

 Since $\log^2\frac{1+8\delta}{1-4\delta} =O(\delta^2)$,
$S^*_{q_3}=O(\delta^2\frac{x}{\log^2x})$.
 Hence for any $\delta$ sufficiently small we get a small number of primes $p\leq
 x$ such that $\frac{p+1}{2}=q_1q_2q_3$ where $q_1$ and $q_2$ are between
$x^{\frac{1}{4}-\delta}$ and $x^{\frac{1}{4}+2\delta}$.
Thus for most such $p$, we have  $q_2>x^{1/4+2\delta}$.

Finally we prove the claim about $q_3$. Assume that
$\frac{p+1}{2}=q_1q_2q_3, \ q_3 \geq q_2 \geq q_1$ then we have
that $q_3 \geq (\frac{p+1}{2})^\frac{1}{3}$.
The following lemma  sharpens this result.
  \begin{lemma}
  For any $0 < \delta < 1/4$ there are at most $O(\delta^2 \frac{x}{\log^2x})$  primes
  $p\leq x$ for which $(\frac{p+1}{2})^\frac{1}{3} \leq q_3 \leq p^{\frac{1}{3}+\delta^2}$
  where $O$ does not depend on $\delta$.
  \end{lemma}
  \begin{proof}
 Note that if $ \frac{p+1}{2}\geq \frac{x}{\log^2x}$ then $q_3 \geq (\frac{p+1}{2})^\frac{1}{3}
   \geq (\frac{x}{\log^2x})^\frac{1}{3}\geq x^{\frac{1}{3}-\delta^2}$ for
 $x\geq x(\delta)$
 (the number of primes $p$ for which
 $\frac{p+1}{2}\leq\frac{x}{\log^2x}$, is $o(\frac{x}{\log^2x})$ by prime number theorem
 and so may be ignored).

 Assume now that $ \frac{p+1}{2}=q_1q_2q_3$ with  $x^{1/3-\delta^2}\leq
 q_3 \leq x^{1/3+\delta^2}$ and $x^{1/4+2\delta}\leq
 q_2 \leq x^{5/12+\delta+\delta^2}$ (this is the maximum range which $q_2$ can be in).
 Using proposition~\ref{proposition 3.5}, we take $a=2q_2q_3, \ b=-1,\
 Y=\frac{x+1}{2q_2q_3}$, and so
\begin{equation*}
\begin{split}
S_{q_1}& =\#\{q_1\leq Y:aq_1+b=\mbox{ prime}\} \\
& =\#\{q_1\leq \frac{x+1}{2q_2q_3}: 2q_2q_3-1=\mbox{ prime}\}\\
& \ll\frac{x+1}{2q_2q_3\log^2\frac{x+1}{2q_2q_3}}
\prod\limits_{\substack{p\mid 2q_2q_3\\p\neq 2}}\frac{p-1}{p-2}
\end{split}
\end{equation*}

Since  $2x^{3/4+\delta+2\delta^2}$ is the maximum of $2q_2q_3$
$(q_1>x^{1/4-\delta})$ we obtain
$$
S_{q_1}\ll
\frac{x}{2q_2q_3\log^2\frac{x}{x^{3/4+\delta}}}\ll
\frac{x}{q_2q_3\log^2x}
$$

Now we sum-up the last term over all possibilities for $q_2,\
q_3$. this number is bounded by, (see in the proof of
lemma~\ref{lemma 3.6})
\begin{equation*}
\begin{split}
S^*_{q_1}& = \frac{x}{\log^2x}
\sum\limits_{x^{\frac{1}{4}+2\delta}
\leq q_2 \leq x^{\frac{5}{12}+\delta+\delta^2}}\frac{1}{q_2}
\sum\limits_{x^{\frac{1}{3}-\delta^2} \leq q_3 \leq
x^{\frac{1}{3}+\delta^2}}\frac{1}{q_3}\\
& \ll
\frac{x}{\log^2x}\log\frac{5/12+\delta
+\delta^2}{1/4+2\delta}\log\frac{1/3+\delta^2}{1/3-\delta^2} \\
& \ll
\frac{x}{\log^2x}\log(1+\frac{6\delta^2}{1-3\delta^2})
= O(\delta^2 \frac{x}{\log^2x})
\end{split}
  \end{equation*}
 and for a sufficiently small $\delta$ we can ignore this number.
 \end{proof}

 By the same method (see lemma 3 in \cite{H}) there are only $O(\delta^2 \frac{x}{\log^2x})$
 primes $p \leq x$ such that $\frac{p+1}{2}=r_1r_2$ where $r_i`s$ are primes,
 $i=1,2,\ \ r_2 \geq r_1, \ \ p^{1/2-\delta^2}\leq r_1 \leq
(\frac{p+1}{2})^{1/2}$.

 If we summarize this section we conclude that for any sufficiently small $0 < \delta < 1/4$ there are at least
 $c_3(\delta)\frac{x}{\log^2x},\ \ c_3(\delta)>0$, primes $p\leq x,\ \ p\equiv u \ (mod
\ v)$ such that we can factor $\frac{p+1}{2}$ in at least one of
the following options:
\begin{enumerate}
\item $\frac{p+1}{2}$ is a prime number.
\item  $\frac{p+1}{2}=r_1r_2$  where $r_1, r_2$ are some
 prime numbers, \ $p^{1/4-\delta}<r_1<p^{1/2-\delta^2},
 \ \ \ p^{1/2+\delta^2}<r_2<p^{3/4+\delta}$.
\item  $\frac{p+1}{2}=q_1q_2q_3$  where  $q_1 \leq q_2 \leq q_3$ are some
  prime numbers, \ $q_1>p^{1/4-\delta}, \ \ \
 q_2>p^{1/4+2\delta},
  \ \ \ q_3 >p^{1/3+\delta^2}$   and so
  $p^{1/2+\delta}<q_1q_2<p^{2/3-\delta^2}, \ \ \ p^{7/12-\delta+\delta^2}<q_1q_3<p^{3/4-2\delta}, \ \ \
  p^{7/12+2\delta+\delta^2}<q_2q_3<p^{3/4+\delta}$.
\end{enumerate}

 \section{Proof of the theorem - The algebraic part}

\subsection{Construction of the arithmetic sequence}

In this section we want to construct integers $u$ and $v$,
$(u,v)=1$ such that for all primes $p$ such that $p\equiv u \ (mod
\ v)$, the discriminants $\Delta_1,...,\Delta_7, \ \Delta_i \neq
3$ of
$\mathbb{Q}(\sqrt{\Delta_1}),...,\mathbb{Q}(\sqrt{\Delta_7})$,
respectively, satisfy
$$
(\frac{\Delta_1}{p})=(\frac{\Delta_2}{p})=...=(\frac{\Delta_7}{p})=-1 \;.
$$
This means that $p$ is inert simultaneously in all of the fields.

In addition we want to insure that $\frac{p+1}{2}$ will be an odd
integer and so we take $u\equiv 1 \ (mod \ 4)$ where $8|v$.
Finally, to get $(\frac{p+1}{2}, v)=1$ we shall construct $u$ and
$v$ so that $(\frac{u+1}{2}, v)=1$ (since after sieving the small
factors of $\frac{p+1}{2}$ we may be left with small factors which
divide $v$, see previous section).

In order to fulfill these demands, we will first show that there
exist infinitely many primes $p$ with the following simultaneous
conditions
\begin{equation}\label{equation 4.1}
 \textstyle(\frac{-1}{p})=(\frac{3}{p})=1 \  \ and \ \
(\frac{\Delta_1}{p})=(\frac{\Delta_2}{p})=...=(\frac{\Delta_7}{p})=-1
\end{equation}

This condition is equivalent to the condition:
$$
B(p)=(1+(\frac{-1}{p}))(1+(\frac{3}{p}))(1-(\frac{\Delta_1}{p}))\cdot
\cdot\cdot (1-(\frac{\Delta_7}{p}))\neq 0
$$

Since the Legendre symbol is a multiplicative function, we obtain,
$$
B(p)=(1+(\frac{-1}{p}))(1+(\frac{3}{p}))(1-\Sigma(\frac{\Delta_i}{p})+
\Sigma(\frac{\Delta_i\Delta_j}{p})-...-(\frac{\Delta_1\cdot\cdot\cdot\Delta_7}{p}))
$$

Let $S$ be the set of all integers of the form
$n=(-1)^{a_1}3^{a_2}\prod\limits_{i=1}^7{\Delta_i^{b_i}}, \ \ a_i,
b_i \in \{0,1\}$. Then
\begin{equation}
\textstyle\sum\limits_{p \leq Z}B(p)=\sum\limits_{n \in
S}(-1)^{b_1+...+b_7}\sum\limits_{p\leq Z}(\frac{n}{p}) \ \ b_i \in
\{0,1\}
\end{equation}

 By the assumption in the theorem (see the introduction) each $n
\in S$ is not a square when $\sum\limits_{i=1}^7b_i$ is odd.

This assumption with the fact that for $n$ not a perfect  square
(by reciprocity law for Legendre symbol)
$$\textstyle\sum\limits_{p\leq Z}(\frac{n}{p})=o(\pi(Z)) \ \ as \
Z \rightarrow \infty$$ implies that $\sum\limits_{p \leq Z}B(p)$
is asymptotic to at least $\pi(Z)$ (since all the negative
summands contribute $o(\pi(Z))$ and at least the natural number 1
contributes $\pi(Z))$. This shows that the simultaneous conditions
have infinitely many solutions $p$.

We fix some particular $p_0$ satisfying the condition. We define
$u_2=p_0$ and for each odd prime $l$, such that
$l|\Delta_1\cdots\Delta_7$ we define $u_l=p_0$ if $l \nmid
p_0+1$ and $u_l=4p_0$ otherwise.
\begin{cl}
$l \nmid u_l+1$
\end{cl}
 \begin{proof}If $u_l=p_0$ then by the assumption $l$
$\nmid p_0+1$, so $l\nmid u_l+1$. If $u_l=4p_0$, assume, by
reductio ad absurdum, that $l|u_l+1$. Hence $l\mid 4p_0+1$.
Because, $u_l=4p_0$ and $l\mid p_0 +1$, we obtain that $l\mid
3p_0$. On the other hand, by our condition, $(\frac{3}{p_0})=1$ so
$(\frac{p_0}{3})=1 \ \ (p_0 \equiv 1 \ (mod \ 4))$. Hence $p_0
\equiv 1 \ (mod \ 3)$. Since $l \mid p_0 +1$ and $p_0 \equiv 1 \
(mod \ 3)$ we conclude that $l \nmid 3$. Using the assumption that
$l \mid p_0 +1$ we deduce that $ l\neq p_0 $ (if $l=p_0$ then $l
\nmid p_0+1$). Hence $l\nmid 3p_0$, a contradiction.
\end{proof}

 Let $v=8\Delta_1\cdots\Delta_7$ and $u$ be the common solution of $u
\equiv u_2 \ (mod \ 8)$ and all the congruences $u\equiv u_l \
(mod \ l)$. Such  a solution exists, by the Chinese Remainder
Theorem.

Since $l\nmid u+1$ for every odd prime $l|v$ and the fact that
$u\equiv 1 \ (mod \ 4)$. (by the construction $u \equiv u_2 \ (mod
 \ 8)$ where $u_2=p_0 \equiv 1 \ (mod \ 4)$) we conclude that $(\frac{u+1}{2},v)=1$.
 Finally, if $p\equiv u \ (mod \ v)$ then, $p\equiv p_0 \ (mod \ 8)$ and $p\equiv p_0$ or $4p_0
 \ (mod \ l)$ for all odd primes $l|v$. So,
$(\frac{\Delta_1}{p})=(\frac{\Delta_1}{p_0})=-1$, and similarly
for all $\Delta_i$`s. This completes the construction of $u$ and
$v$.

Note that by the construction of the integers $u$ and $v$ we have
that $(u,v)=1$. (take $l$ an odd prime number, $l\mid
v=8\Delta_1\cdots\Delta_7$ and assume that $l\mid u$. Since $u
\equiv u_l \ (mod \ l)$ then $l \mid u_l$ hence $l \mid p_0$ or
$4p_0$ in other words $l=p_0$. But $p_0 \nmid
\Delta_1\cdots\Delta_7$ ($p_0$ fulfills the simultaneous condition
\eqref{equation 4.1}) and $l \mid \Delta_1\cdots\Delta_7$).

\subsection{The last step of the proof}

For the last step of the proof we need to use lemma 4 from
Narkiewicz  \cite{N1}, which generalized lemma 2 in \cite{GM}.
 \begin{lemma}\label{lemma 4.2}
If $a_1,\ldots a_k$ are multiplicatively independent integers of an
 algebraic number-field $K$, $G$ the subgroup of $K^\star$
 generated by $a_1,\ldots a_k$, and for any prime ideal
 $\mathbf{P}$ not dividing $a_1,\cdots a_k$ we denote by
 $G_{\mathbf{P}}$ the reduction of $G \ (mod  \ \mathbf{P})$, then for
 all positive $y$ one can have $\#G_{\mathbf{P}} < y $ for at most
 $O(y^{1+\frac{1}{k}})$ prime ideals $\mathbf{P}$, with the implied
 constant being dependent on the $a_i$`s and $K$.
 \end{lemma}

Now, as we saw at the end of section 3, for any sufficiently small
$0 <\delta <1/4$ there is some constant $c_3(\delta) > 0$ such
that for $c_3(\delta) \frac{x}{\log^2x}$ primes $p \leq x, \ p
\equiv u \ (mod \ v)$ at least one of the following occur:
\begin{enumerate}
\item $\frac{p+1}{2}$ is a prime number.
\item $\frac{p+1}{2}=r_1r_2$  where $r_1, r_2$ are
 primes so that,  \ $p^{1/4-\delta}<r_1<p^{1/2-\delta^2},
 \ \ \ p^{1/2+\delta^2}<r_2<p^{3/4+\delta}$.
\item $\frac{p+1}{2}=q_1q_2q_3$ where $q_1 \leq q_2 \leq q_3$  are
  primes such that, $p^{1/2+\delta}<q_1q_2<p^{2/3-\delta^2}$,
$p^{7/12-\delta+\delta^2}<q_1q_3<p^{3/4-2\delta}$,
  $p^{7/12+2\delta+\delta^2}<q_2q_3<p^{3/4+\delta}$.
\end{enumerate}
It is clear by the construction of $u$ and $v$ that $p\equiv 1 \
(mod \ 4)$. Because $\#C_\epsilon(p)=p+1$ when $p$ is inert in
${\mathbb{Q}(\sqrt{\Delta})}$ the unit $-1$ is a non-square in the
group $C_\epsilon(p)$. Hence for any unit $\epsilon$ , we can
choose constant $c=\pm1$ such that $c\epsilon$ is a non-square in
$C_\epsilon(p)$. Similarly, since $c\epsilon$ is a non-square and
the index of the group of squares is 2, by the theorem
on cyclic groups, the order of $c\epsilon$ is even.

Now we look at our cases:\\

(1) In this case, by the above note, $c\epsilon$, if not
primitive, has order 2 But the number of $p$`s with this property
is O(1) (by lemma~\ref{lemma 4.2}).\\

(2) Let $c_1\epsilon_1,\dots ,c_4\epsilon_4$  be units in the orders
 $\mathcal{O}_{\Delta_1}$, $\dots ,\mathcal{O}_{\Delta_4}$  of
${\mathbb{Q}(\sqrt{\Delta_1})}$, $\dots, {\mathbb{Q}(\sqrt{\Delta_4})}$,
 $\Delta_i\neq 3$,  $i=1,2,3,4$,
respectively. We will show that one of them is primitive
infinitely many times.

If $\ord(c_i\epsilon_i)=2r_1 < 2x^{1/2-\delta^2}$ for some
$i=1,2,3,4$. by Lemma~\ref{lemma 4.2} this occurs in at most
$O(x^{1/2-\delta^2})^2=x^{1-2\delta^2}$ and this is a negligible
number compared to $c_3(\delta) \frac{x}{\log^2x}$ .

Assume $\ord(c_i\epsilon_i)=2r_2, \ i=1,2,3,4$. Consider the ring
of integers $\mathcal{O}_M$ of the compositum field $M$ of
$\mathbb{Q}(\sqrt{\Delta_i}) \ i=1,...,4$.

\begin{prop}
For any prime ideal $P \mid (p)=p\mathcal{O}_M$:
$$
(\mathcal{O}_M / \mathbf{P})^\star\simeq (\mathcal{O}_{\Delta_1}/
p\mathcal{O}_{\Delta_1})^\star\simeq ...\simeq
(\mathcal{O}_{\Delta_4}/ p\mathcal{O}_{\Delta_4})^\star \;.
$$
\end{prop}

\begin{proof}
Since $p$ is inert, the order of
$\mathcal{O}_{\Delta_i}/p\mathcal{O}_{\Delta_i}$ is $p^2$, i.e.,
$[\mathcal{O}_{\Delta_i}/p\mathcal{O}_{\Delta_i}:\mathbb{Z}/p\mathbb{Z}]=2.$
Since all these quotient fields are finite fields and two finite
fields with the same number of elements are isomorphic, it is
enough to show that,
$$ f= [\mathcal{O}_{M}/P:\mathbb{Z}/p\mathbb{Z}]=2. $$

Consider the Galois group \ $G=Gal[M/\mathbb{Q}]$ \ and define two
subgroups of $G$, the decomposition group $D$ and the inertia
group $E$:
$$ D=D(P|(p))=\{\sigma \in G|\ \sigma(P)=P \}$$
and
$$ E=E(P|(p))=\{\sigma \in G|\ \sigma(\alpha) \equiv \alpha
\ (mod \ P),\  \forall \alpha \in \mathbb{Z} \}\;.
$$

Now, consider the Galois group $\bar{G}$,
$$ \bar{G}=Gal[\mathcal{O}_{M}/P \ / \ \mathbb{Z}/p\mathbb{Z}]$$

By  \cite[ chapter 4, beginning]{Mar},
$$ D/E \simeq \bar{G} $$
By theorem 28 in in \cite{Mar} (since $(p)$ is inert in all the
fields $\mathbb{Q}(\sqrt{\Delta_i})$, $(p)$ is unramified in all
the fields $\mathbb{Q}(\sqrt{\Delta_i})$. Hence $(p)$ is also
unramified in $\mathcal{O}_M$, i.e., $e$ - the exponent of $P$ in
the decomposition of $(p)$, equal to 1)
$$|E|=e=1$$
$$|D|=f$$
Immediately we conclude that,
$$ D \simeq \bar{G}$$
Since $\bar{G}$ is a cyclic group, we get that $D$ is a cyclic
group of order $f$. But $D$ is a subgroup of $G$ and $G=C_2 \times
...\times C_2$ where $C_2$ is a group of order 2. So, $f\leq 2$.
Since
$$
[\mathcal{O}_{\Delta_i}/p\mathcal{O}_{\Delta_i}:\mathbb{Z}/p\mathbb{Z}]=2
$$
we also see that $f\geq 2$. Hence $f=2$.

Because the compositum of normal extensions is normal, this claim
is true for all the prime ideals $P$ in the decomposition of $(p)$
(they have the same $e$ and the same $f$)
\end{proof}

 By the last proposition  $|<c_1\epsilon_1,c_2\epsilon_2,c_3\epsilon_3,c_4\epsilon_4>|=2r_2
  \ll p^{3/4+\delta}$ in $(\mathcal{O}_M/ \mathbf{P})^\star$  for $(p)|P$
  (this is a cyclic group). By Lemma~\ref{lemma 4.2}
the number of $p \leq x$ that have this order is
at most $O(x^{3/4+\delta})^\frac{5}{4}=O(x^{15/16+5\delta/4})$ but
we can choose  $\delta$ to be as small as needed. Hence the set of
primes $p,\ p \leq x$ such that $c_i\epsilon_i \ i=1,...,4$ have
order $2r_2$ is a small in comparison with $c_3(\delta)
  \frac{x}{\log^2x}$. \\

(3) As we did in part (2), we define $c_i\epsilon_i, \ c_i=\pm
1$ in the ring of integers $\mathcal{O}_{\Delta_i}$ of
${\mathbb{Q}(\sqrt{\Delta_i})},\ \Delta_i\neq3, \ i=1,...,7$ ,
respectively. Assume that each has order $2q_j < 2x^{1/2-\delta},
\ j=1,2$ or $3$. By lemma~\ref{lemma 4.2} at most $O(x^{1-2\delta})$ primes
$p\leq x$ have this property. Now, let`s take any two units of the
seven units above, $c_j\epsilon_j$ and $\ c_k\epsilon_k \ \ 1 \leq
j,k \leq 7 $ and assume that they have order $2q_1q_2 <
2x^{2/3-\delta^2}$. By lemma~\ref{lemma 4.2}, we can prove (in the same way we
did in (2)), that this occurs for at most $O(x^{1-3/2\delta^2})$
primes $p \leq
x$.\\

To summarize, we take $\delta > 0$ such that for $c_3(\delta)
\frac{x}{\log^2x}$ primes $p \leq x$,
$c_1\epsilon_1,...,c_6\epsilon_6$ do not have order $2q_j, \
j=1,2,3$ and $2q_1q_2$.

 From these six units we take three,
$c_{i_1}\epsilon_{i_1},c_{i_2}\epsilon_{i_2},c_{i_3}\epsilon_{i_3}$
and assume that they have order $2q_1q_3 < 2x^{3/4-2\delta}$.
Again by lemma~\ref{lemma 4.2} this occurs for at most $O(x^{1-8/3\delta})$
primes $p\leq x$. In other words at least four of them do not have
order $2q_1q_3$ for $c_3(\delta) \frac{x}{\log^2x}$ primes.

Consider these four units $c_1\epsilon_1,...,c_4\epsilon_4$
(without loss generality) and assume that they have an order
$2q_2q_3 < 2x^{3/4+\delta}$. Again by lemma~\ref{lemma 4.2} this
occur in at most $O(x^{15/16+5/4\delta})$ primes $p\leq x$. So,
for a sufficiently small $0 < \delta < 1/4$, there is at least one
unit, say $c_1\epsilon_1$, such that for $c_3(\delta)
\frac{x}{\log^2x}$ primes $p \leq x$, $c_1\epsilon_1$ is
primitive.\\

Note that  this theorem implies Corollary~\ref{corollary 1.7}. \\

{\bf Acknowledgement} I wish to gratefully acknowledge my Ph.D
supervisors Prof. Zee`v Rudnick and Prof. Jack Sonn for the
helpful suggestion and the fruitful ideas.

\end{document}